\documentclass{article}
\usepackage{amssymb,latexsym}
\usepackage{amsmath}
\usepackage{graphics}
\usepackage{graphicx}
\newtheorem{theorem}{Theorem}

\begin{document}
\title{ONE-PARAMETER HOMOTHETIC MOTION IN THE HYPERBOLIC PLANE AND EULER-SAVARY FORMULA }
\author{Soley ERSOY, Mahmut AKYIGIT\\
sersoy@sakarya.edu.tr , makyigit@sakarya.edu.tr\\
         Department of Mathematics,
Faculty of Arts and Sciences
\\
Sakarya University, 54187 Sakarya/TURKEY \\}

\date{March 17, 2009}

\maketitle

\begin{abstract}
In \cite{Mul} one-parameter planar motion was first introduced and
the relations between absolute, relative, sliding velocities (and
accelerations) in the Euclidean plane  $\mathbb{E}^2$ were
obtained. Moreover, the relations between the Complex velocities
one-parameter motion in the Complex plane were provided by
\cite{Mul}. One-parameter planar homothetic motion was defined in
the Complex plane, \cite{Kur}. In this paper, analogous to
homothetic motion in the Complex plane given by \cite{Kur},
one-parameter planar homothetic motion is defined in the
Hyperbolic plane. Some characteristic properties about the
velocity vectors, the acceleration vectors and the pole curves are
given. Moreover, in the case of homothetic scale $h$ identically
equal to 1, the results given in \cite{Yuc} are obtained as a
special case. In addition, three hyperbolic planes, of which two
are moving and the other one is fixed, are taken into
consideration and a canonical relative system for one-parameter
planar hyperbolic homothetic motion is defined. Euler-Savary
formula, which gives the relationship between the curvatures of
trajectory curves, is obtained with the help of this relative
system.


\end{abstract}

\section{Preliminaries}\label{S:intro}

Before proceeding any further, we require a definition for the set
of hyperbolic number and assume the existence of any number $j$
which has the property $ j \ne  \pm 1$ . In terms of the standard
basis $\left\{ {1,j} \right\}$, the hyperbolic number  can be
written as
\[
z = x + jy
\]
where $j\left( {j^2  = 1} \right)$ is the unipotent (hyperbolic)
imaginary unit and the reel numbers $x$ and $y$ are called the
real and unipotent (or hallucinatory) parts of the hyperbolic
number $z$, respectively,
\cite{Cat}-\cite{Coc},\cite{Fje}-\cite{Fje3},\cite{Roc},
\cite{Sob}. The hyperbolic numbers
\[
 \mathbb{H}= \mathbb{R} \left[ j \right] = \left\{ {\left. {z = x + jy} \right|x,y \in ,j^2  = 1} \right\}
\]
are the real numbers extended to include the unipotent $j$  in the
same manner that $ \mathbb{C}= \mathbb{R}\left[ i \right]$ are the
complex numbers extended to include the imaginary $i$, $\left(
{i^2  =  - 1} \right)$, \cite{Yag}.\\
The hyperbolic numbers are also called perplex numbers \cite{Fje},
split-complex numbers \cite{Alf} or double numbers
\cite{Alf},\cite{Fje2},\cite{Roc}. The hyperbolic number systems
serve as the coordinates in the Lorentzian plane in the same way
as the complex numbers serve as coordinates in the Euclidean
plane. The role played by the complex numbers in Euclidean space
is played by the hyperbolic number systems in the pseudo-Euclidean
space, \cite{Sob}.\\
Addition and multiplication of the hyperbolic numbers are
\[
\displaylines{
  \left( {x + jy} \right) + \left( {u + jv} \right) = \left( {x + u} \right) + j\left( {y + v} \right), \cr
  \left( {x + jy} \right)\left( {u + jv} \right) = \left( {xu + yv} \right) + j\left( {xv + yu} \right). \cr}
\]
respectively. This multiplication is commutative, associative and
distributes over addition. The hyperbolic conjugate of $z = x +
jy$ is defined by $\overline z  = x - jy$. The hyperbolic inner
product is
\[
\left\langle {z,w} \right\rangle  = {\mathop{\rm Re}\nolimits}
\left( {z\overline w } \right) = {\mathop{\rm Re}\nolimits} \left(
{\overline z w} \right) = xu - yv
\]
where; $z = x + jy$ and $w = u + jv$. Hyperbolic numbers $z$ and
$w$ are hyperbolic (Lorentzian) orthogonal if  $ \left\langle
{z,w} \right\rangle  = 0$, \cite{Sob}.\\
The hyperbolic modulus of $z = x + jy$ is
\[
\left\| z \right\|_h  = \sqrt {\left| {\left\langle {z,z}
\right\rangle } \right|}  = \sqrt {\left| {z\overline z } \right|}
= \sqrt {\left| {x^2  - y^2 } \right|}
\]
and it is the hyperbolic distance of the point $z$ from the
origin. This is the Lorentz invariant of two-dimensional special
relativity and their unimodular multiplicative group (the group
composed of quadratic matrices determinant of which equals to 1)
is the special relativity Lorentz group, \cite{Yag2}. These
relations have been used to extend special relativity.
Furthermore, by using the functions of the hyperbolic variable,
two-dimensional special relativity has been generalized,
\cite{Cat2}. These applications make the hyperbolic numbers
appropriate for physics and the application of hyperbolic numbers
is similar to the application of complex numbers to the Euclidean
plane geometry, \cite{Yag2}.\\
Note that the points $z \ne 0$ on the lines $y = x$ are isotropic
in the sense that they are nonzero vectors with $ \left\| z
\right\|_h  = 0$. By this way, the hyperbolic distance creates
Lorentzian geometry in $\mathbb{R}^2$. This is different from the
usual Euclidean geometry of the complex plane, where $\left\| z
\right\|_h  = 0$ only if $z=0$ in the complex plane. The set of
all points in the hyperbolic plane that satisfy the equation
$\left\| z \right\|_h  = r > 0$ is a four-branched hyperbola of hyperbolic radius $r$, \cite{Sob}.\\
The hyperbolic number $z = x + jy$ can be written as follows:\\
While the hyperbolic number $z$ is on H-I or H-III plane, then

$$z = \pm r\left( {\cosh \varphi  + j\sinh \varphi } \right) =  \pm
re^{j\varphi },$$ While the hyperbolic number $z$ is on H-II or
H-IV plane, then
$$z = \pm r\left( {\sinh \varphi  + j\cosh \varphi
} \right) =  \pm rje^{j\varphi },$$ [See Figure 1.]

\begin{center}
\hfil\scalebox{1}{\includegraphics{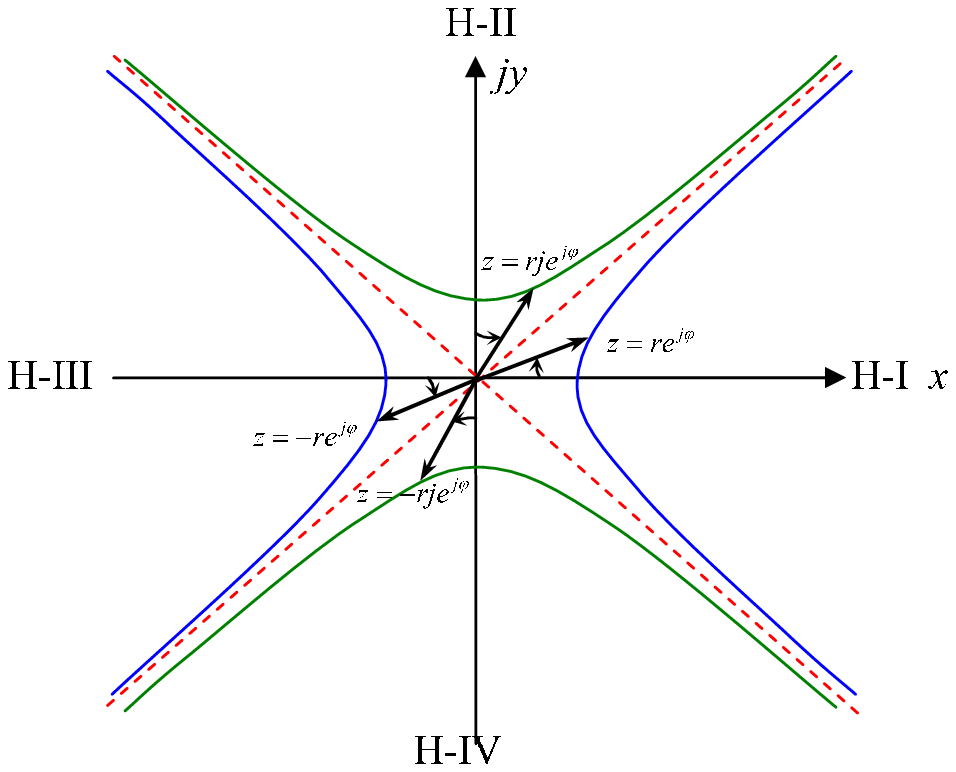}}\hfil\\
\scriptsize{Figure 1. Hyperbolic Plane}\\
\end{center}
\normalsize This formula can be derived by a power
series expansion due to the fact that cosh has only even powers
whereas sinh has odd powers. For all real values of the hyperbolic
angle $\varphi$, the hyperbolic number $e^{j\varphi }$ has norm 1
and lies on the right
branch of the unit hyperbola, \cite{Sob}.\\
A hyperbolic rotation defined by $e^{j\varphi }$ corresponds to
multiplication by the matrix, \cite{Sob};
\[
\left[ {\begin{array}{*{20}c}
   {\cosh \varphi } & {\sinh \varphi }  \\
   {\sinh \varphi } & {\cosh \varphi }  \\
\end{array}} \right].
\]
Another property of the hyperbolic inner product is
\[
\left\langle {ze^{j\varphi } ,we^{j\varphi } } \right\rangle  =
\left\langle {z,w} \right\rangle .
\]
In addition, a vector multiplied by $j$ is a hyperbolic orthogonal
vector. This is similar to the role played by the multiplication $
i = e^{i\left( {{\pi  \mathord{\left/
 {\vphantom {\pi  2}} \right.
 \kern-\nulldelimiterspace} 2}} \right)}$ in the complex plane, \cite{Sob}.

\section{One-Parameter Homothetic Motion in the Hyperbolic Plane }
In this section, we will define one-parameter homothetic motion in
the hyperbolic plane and obtain the relation between the
velocities and the accelerations of a point under one-parameter
homothetic planar motions.\\
Homothetic motion of a moving hyperbolic plane $\mathbb{H}$  with
respect to a fixed hyperbolic plane $\mathbb{H'}$ will be
considered, in that, the orthonormal coordinate systems $\left\{
{O;{\bf h}_{\bf 1} ,{\bf h}_{\bf 2} } \right\}$ and $\left\{
{O';{\bf h'}_1 ,{\bf h'}_2 } \right\}$ being on the moving and
fixed hyperbolic planes $\mathbb{H}$ and $\mathbb{H'}$,
respectively, will be analyzed with respect to each other. As the
vector $ \overrightarrow {OO'}$ represented by the hyperbolic
number $ {\bf u}$ determines the distance between the origin point
of the moving system and the origin point of the fixed system, the
vectorial representation is as follows
\begin{equation}\label{2.1}
{\bf x'} = h{\bf x} - {\bf u}
\end{equation}
where, $h = h\left( t \right) \ne {\rm constant} $ is the
homothetic scale (Figure 2. and Figure 3.).\\
\hfil\scalebox{1}{\includegraphics{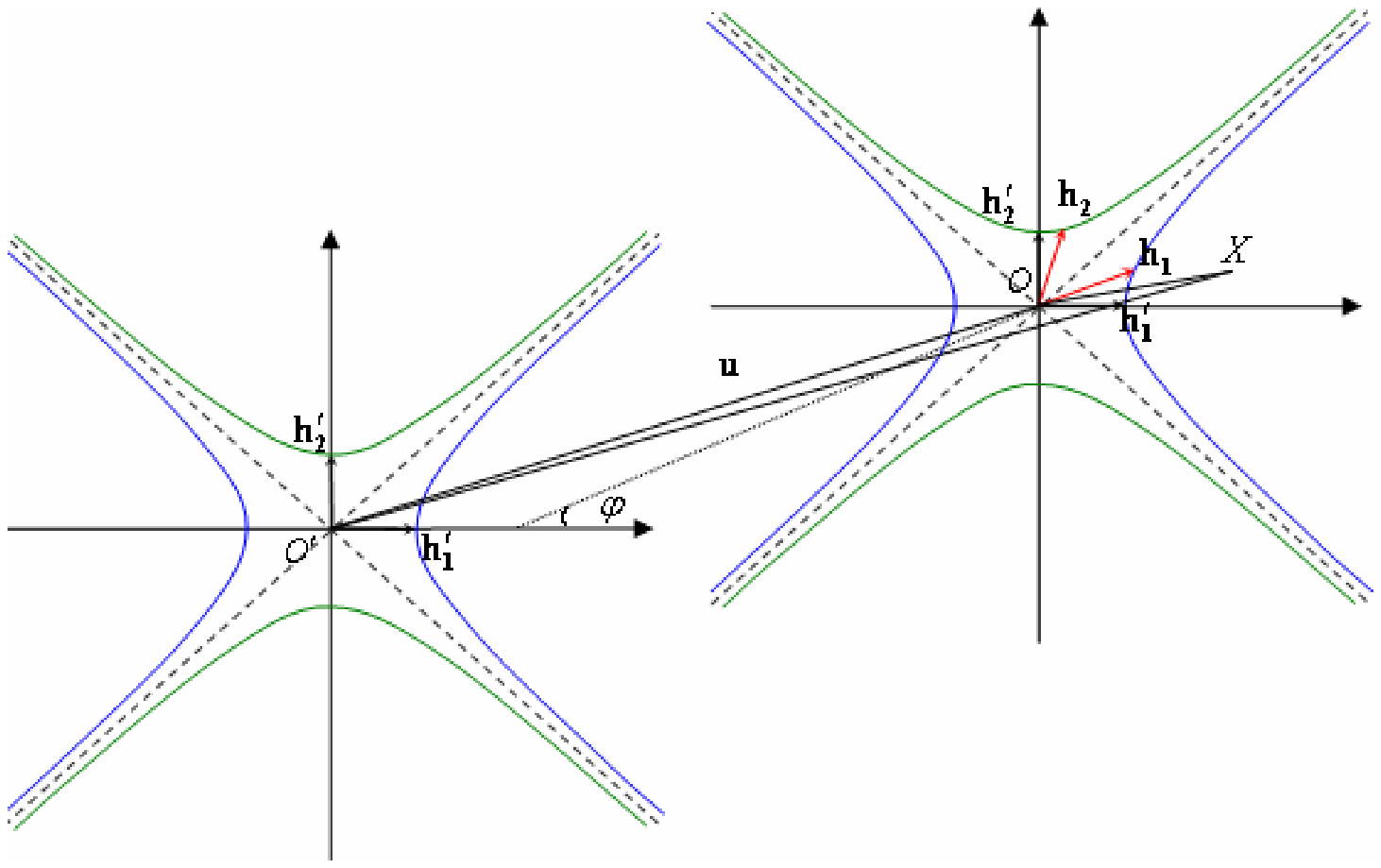}}\hfil
\begin{center}
\scriptsize{Figure 2. $ {\mathbb{H} \mathord{\left/
 {\vphantom {a b}} \right.
 \kern-\nulldelimiterspace} \mathbb{H'}}$ hyperbolic homothetic motion that rotates with central angle $\varphi$}
\end{center}
\hfil\scalebox{1}{\includegraphics{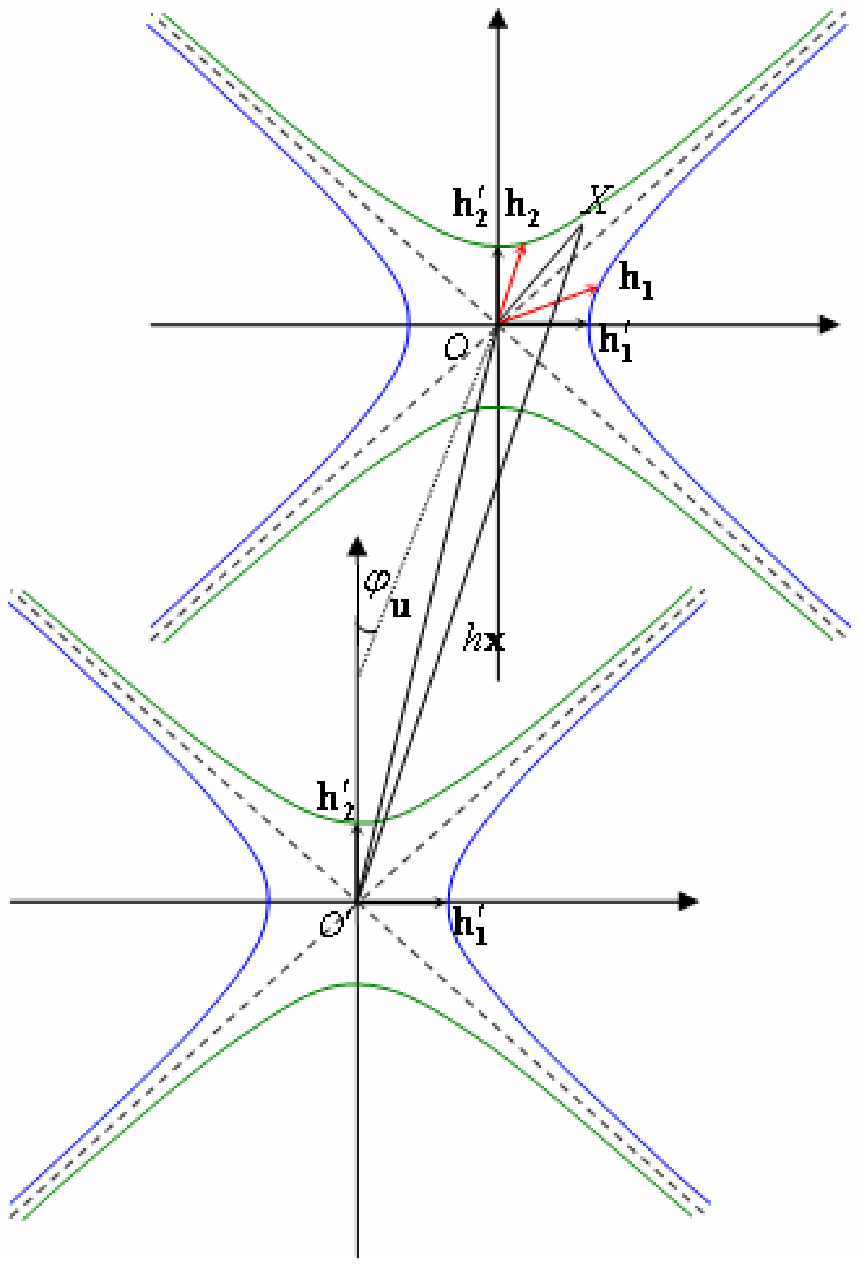}}\hfil
\begin{center}
\scriptsize{Figure 3.  $ {\mathbb{H} \mathord{\left/
 {\vphantom {a b}} \right.
 \kern-\nulldelimiterspace} \mathbb{H'}}$ hyperbolic homothetic motion that rotates with hyperbolic angle  $\varphi$}
\end{center}
\normalsize
Let a fixed point $X$  chosen on the plane
$\mathbb{H}$ be represented by the hyperbolic numbers ${\bf x} =
x_1  + jx_2$ and ${\bf x'} = x'_1  + jx'_2$ on the planes
$\mathbb{H}$ and $\mathbb{H'}$, respectively. Thus, one-parameter
homothetic hyperbolic motion of the moving coordinate system
$\left\{ {O;{\bf h}_{\bf 1} ,{\bf h}_{\bf 2} } \right\}$  with
respect to the fixed coordinate system $\left\{ {O';{\bf h'}_1
,{\bf h'}_2 } \right\}$ represented by $ {\mathbb{H}
\mathord{\left/
 {\vphantom {a b}} \right.
 \kern-\nulldelimiterspace} \mathbb{H'}}$ is defined as the following transformation;
\begin{equation}\label{2.2}
{\bf x'} = {\bf u'} + h{\bf x}e^{j\varphi } .
\end{equation}
Here, $\varphi$ is the Lorentzian (either central or hyperbolic)
rotation angle of the motion $ {\mathbb{H} \mathord{\left/
 {\vphantom {a b}} \right.
 \kern-\nulldelimiterspace} \mathbb{H'}}$, and the hyperbolic
 number $ {\bf u'}$ represents the origin point of the moving system in the fixed
system. The rotation angle $\varphi$, the homothetic scale $h$ and
${\bf x}{\bf ,}\;{\bf x'}{\bf ,}\;{\bf u}$ will be regarded as the
differentiable functions of a real parameter $t$ from the class $
C^\infty $. Generally, this parameter $t$ will be used as time and
at the moment $t=0$, the coordinate systems will be accepted as
coincident.\\ The hyperbolic number ${\bf u} = u_1  + ju_2$
represents the origin point $O'$ of the fixed system in the moving
system. At this case, if $X'=O'$, then $ {\bf x' = 0}$ and ${\bf
x=u}$. Thus, from the equation (\ref{2.2})
\begin{equation}\label{2.3}
{\bf u'} =  - {\bf u}e^{j\varphi }
\end{equation}
is found. Using (\ref{2.2}) and (\ref{2.3}) the following is
obtained.
\begin{equation}\label{2.4}
{\bf x'} = \left( {h{\bf x} - {\bf u}} \right)e^{j\varphi } .
\end{equation}
As $\dot \varphi \left( t \right)$ would give only the
translation, we will assume $\frac{{d\varphi }}{{dt}} = \dot
\varphi \left( t \right) \ne 0 $ during the motion $ {\mathbb{H}
\mathord{\left/
 {\vphantom {a b}} \right.
 \kern-\nulldelimiterspace} \mathbb{H'}}$ and call it
as the angular velocity of the motion.
\subsection{Velocities and the Composition of Velocities}
Let the point $X$ on the moving plane $\mathbb{H}$ change its
location depending on a parameter $t$ while undergoing
one-parameter homothetic motion of plane $\mathbb{H}$ with regards
to the plane $\mathbb{H'}$. At this case, two motions belonging to
the point $X$ occur. Let's explore what kind of relation exists
between these homothetic motions and their velocities.\\
The velocity vector of the point $X$ with respect to the plane
$\mathbb{H}$, that is, the vectorial velocity which the point has
while drawing the trajectory curve on is called the relative
velocity of the point and is represented as ${\bf V}_{\bf r}$. The
relative velocity ${\bf V}_{\bf r}$ of $X$ is
\begin{equation}\label{2.5}
{\bf V}_{\bf r}  = h{\bf \dot x}e^{j\varphi } .
\end{equation}
The velocity of the point $X$ with respect to the fixed plane
$\mathbb{H'}$ is called the absolute velocity of the point $X$ and
is represented by ${\bf V}_{\bf a}$. The absolute velocity of $X$
is
\begin{equation}\label{2.6}
{\bf V}_{\bf a}  = \left( {\dot h + jh\dot \varphi } \right){\bf
x}e^{j\varphi }  - \left( {{\bf \dot u} + j{\bf u}\dot \varphi }
\right)e^{j\varphi }  + h{\bf \dot x}e^{j\varphi }
\end{equation}
If the differential of equation (\ref{2.3}) is substituted in the
last equation
\begin{equation}\label{2.7}
{\bf V}_{\bf a}  = {\bf \dot u'} + \left( {\dot h + jh\dot \varphi
} \right){\bf x}e^{j\varphi }  + h{\bf \dot x}e^{j\varphi }
\end{equation}
can be written. If
\begin{equation}\label{2.8}
{\bf V}_{\bf f}  = \left( {\dot h + jh\dot \varphi } \right){\bf
x}e^{j\varphi }  - \left( {{\bf \dot u} + j{\bf u}\dot \varphi }
\right)e^{j\varphi }
\end{equation}
the following equation is obtained
\begin{equation}\label{2.9}
{\bf V}_{\bf a}  = {\bf V}_{\bf f}  + {\bf V}_{\bf r} .
\end{equation}
Here, ${\bf V}_{\bf f}$ is the sliding velocity of the
one-parameter planar homothetic motion $ {\mathbb{H}
\mathord{\left/
 {\vphantom {a b}} \right.
 \kern-\nulldelimiterspace} \mathbb{H'}}$. If the point $X$ is a fixed
point on the moving plane $\mathbb{H}$, then  ${\bf V}_{\bf
f}={\bf 0}$. So it is easily seen that
\begin{equation}\label{2.10}
{\bf V}_{\bf a}  = {\bf V}_{\bf r} .
\end{equation}
\subsection{The Rotation Pole and the Pole Trajectories}
Studying the points of the one-parameter homothetic hyperbolic
motion $ {\mathbb{H} \mathord{\left/
 {\vphantom {a b}} \right.
 \kern-\nulldelimiterspace} \mathbb{H'}}$ where the ${\bf V}_{\bf
f}$ sliding velocity equals to zero in every moment $t$ will
reveal the rotation pole term. Thus, by taking ${\bf V}_{\bf
f}={\bf 0}$ in the equation (\ref{2.8}), the pole point $P =
\left( {p_1 ,p_2 } \right) \in \mathbb{H}$  is the hyperbolic
number as
\begin{equation}\label{2.11}
{\bf p} = \frac{{{\bf \dot u} + j\dot \varphi {\bf u}}}{{\dot h +
jh\dot \varphi }}
\end{equation}
or
\[
{\bf p} = p_1  + jp_2  = \frac{{\dot h{\bf \dot u} - h\dot \varphi
^2 {\bf u}}}{{\dot h^2  - h^2 \dot \varphi ^2 }} + j\frac{{\dot
h\dot \varphi {\bf u} - h\dot \varphi {\bf \dot u}}}{{\dot h^2  -
h^2 \dot \varphi ^2 }}.
\]
In the special case of $h\left( t \right) = 1$, the following
equation exists:
\[
{\bf p} = p_1  + jp_2  = {\bf u} + j\frac{{{\bf \dot u}}}{{\dot
\varphi }}
\]
which was given in \cite{Yuc}.\\
Let the rotation pole of the homothetic hyperbolic motion $
{\mathbb{H} \mathord{\left/
 {\vphantom {a b}} \right.
 \kern-\nulldelimiterspace} \mathbb{H'}}$ be $P$
and a moving point on $\mathbb{H'}$ be $X$. Given this condition,
the pole ray $ \overrightarrow {PX}$ from the pole $P$ to the
point $X$ is expressed by the following equation
\begin{equation}\label{2.12}
\overrightarrow {PX}  = \left( {h{\bf x} - {\bf p}}
\right)e^{j\varphi } .
\end{equation}
In addition, if the equations (\ref{2.8}) and (\ref{2.11}) are
considered together
\begin{equation}\label{2.13}
{\bf V}_{\bf f}  = \left( {\dot h + jh\dot \varphi } \right)\left(
{{\bf x - p}} \right)e^{j\varphi }
\end{equation}
is found. As in \cite{Yuc}, ${\bf V}_{\bf f}$ and $
\overrightarrow {PX}$ are seen to be orthogonal to each other
given condition of $h\left( t \right) =
1$.\\
The length of the vector obtained from equation (\ref{2.13}) is
\[
\left\| {{\bf V}_{\bf f} } \right\|_h  = \sqrt {\left( {\dot h^2
- h^2 \dot \varphi ^2 } \right)\left[ {\left( {x_1  - p_1 }
\right)^2  + \left( {x_1  - p_1 } \right)^2 } \right]} .
\]
In special case of  $h\left( t \right) = 1$ as given in
\cite{Yuc},
\[
\left\| {{\bf V}_{\bf f} } \right\|_h  = \left| {\dot \varphi }
\right|\left\| {\overrightarrow {PX} } \right\|_h
\]
is found.\\
During one-parameter homothetic hyperbolic motion $ {\mathbb{H}
\mathord{\left/
 {\vphantom {a b}} \right.
 \kern-\nulldelimiterspace} \mathbb{H'}}$ the geometric
locus of the pole points $P$ in each $t$ moment is the moving pole
curve $ \left( P \right)$ on the plane $ \mathbb{H}$ and the fixed
pole curve $ \left( P' \right)$ on the plane $ \mathbb{H'}$,
respectively. Due to equation (\ref{2.1}), the following is
written:
\[
{\bf p'} = \left( {h{\bf p} - {\bf u}} \right)e^{j\varphi }
\]
and from the differentiation of this last equation with respect to
$t$
\[
{\bf \dot p'} = \left[ {\left( {\dot h + jh\dot \varphi }
\right){\bf p} - \left( {{\bf \dot u} + j{\bf u}\dot \varphi }
\right) + h{\bf \dot p}} \right]e^{j\varphi }
\]
is obtained. Here, if the equation of the pole point given by
equation (\ref{2.11}) is substituted in the last equation,
\begin{equation}\label{2.14}
{\bf \dot p'} = h{\bf \dot p}e^{j\varphi }
\end{equation}
is found.\\
Thus, the tangent vectors at the contact points of the pole curves
coincide with each other after the Lorentzian rotation $\varphi$
and the translation $h$.\\
Let the arc elements of the moving and the fixed pole curves be
$ds$ and $ds'$, respectively. In this case
\[
ds = \left\| {{\bf \dot p}} \right\|_h dt\quad {\rm and}\quad ds'
= \left\| {{\bf \dot p'}} \right\|_h dt
\]
can be written. With the help of this last equation and equation
(\ref{2.14}), we get
\[
ds' = \left| h \right|ds.
\]
Thus, the following theorem can be given.
\begin{theorem}\label{T:2.1}
In one-parameter planar homothetic hyperbolic motion $ {\mathbb{H}
\mathord{\left/
 {\vphantom {a b}} \right.
 \kern-\nulldelimiterspace} \mathbb{H'}}$, the moving
pole curve $ \left( P \right)$ in the plane $\mathbb{H}$ rolls by
sliding on the fixed pole curve $ \left( P' \right)$ on the plane
$\mathbb{H'}$. The coefficient of this sliding, rolling motion is
the homothetic scale $h$.
\end{theorem}
\textbf{Special Case} In the special case of $h=1$, we get
$ds'=ds$, that is, the pole curves roll on each other without
sliding, which is given in \cite{Yuc}.
\subsection{Accelerations and the Composition of Accelerations}
Let $X$ be a moving point on the moving hyperbolic plane
$\mathbb{H}$. In this case, the acceleration of the point $X$ with
respect to $\mathbb{H}$ is called relative acceleration and is
defined by $ \frac{{d^2 {\bf x}}}{{dt^2 }} = {\bf \ddot x}$. In
addition, the relative acceleration ${\bf b}_{\bf r}$ with respect
to fixed hyperbolic plane $\mathbb{H'}$ can be written as
\begin{equation}\label{2.15}
{\bf b}_{\bf r}  = h{\bf \ddot x}e^{j\varphi } .
\end{equation}
The acceleration of $X'$ with respect to $\mathbb{H'}$ is called
the absolute acceleration ${\bf b}_{\bf a}$, and with the help of
the differentiation of ${\bf V}_{\bf a}$ with respect to $t$,
\begin{equation}\label{2.16}
{\bf b}_{\bf a}  = \frac{{d{\bf V}_{\bf a} }}{{dt}} = {\bf \dot
V}_{\bf a}  = \left( {{\bf x - p}} \right)\left[ {\ddot h + h\dot
\varphi ^2  + j\left( {2\dot h\dot \varphi  + h\ddot \varphi }
\right)} \right]e^{j\varphi }  - {\bf \dot p}\left( {\dot h +
jh\dot \varphi } \right)e^{j\varphi }  + 2{\bf \dot x}\left( {\dot
h + jh\dot \varphi } \right)e^{j\varphi }  + h{\bf \ddot
x}e^{j\varphi }
\end{equation}
is obtained and in this last equation, the
expression
\begin{equation}\label{2.17}
{\bf b}_{\bf f}  = \left( {{\bf x - p}} \right)\left[ {\ddot h +
h\dot \varphi ^2  + j\left( {2\dot h\dot \varphi  + h\ddot \varphi
} \right)} \right]e^{j\varphi }  - {\bf \dot p}\left( {\dot h +
jh\dot \varphi } \right)e^{j\varphi }
\end{equation}
is called the sliding acceleration and
\begin{equation}\label{2.18}
{\bf b}_{\bf c}  = 2{\bf \dot x}\left( {\dot h + jh\dot \varphi }
\right)e^{j\varphi }
\end{equation}
is called the Coriolis acceleration.\\
Thus, the following theorem can be given.
\begin{theorem}\label{T:2.2}
In $ {\mathbb{H} \mathord{\left/
 {\vphantom {a b}} \right.
 \kern-\nulldelimiterspace} \mathbb{H'}}$ one-parameter planar homothetic hyperbolic motion, there is
the following relation between the accelerations
\[
{\bf b}_{\bf a} {\bf  = b}_{\bf f} {\bf  + b}_{\bf c} {\bf  +
b}_{\bf r}.
\]
\end{theorem}
The acceleration pole is known by the vanishing of the sliding
acceleration under one-parameter planar motion. Thus, the
following theorem is obtained, given the condition of ${\bf
b}_{\bf f}={\bf 0}$.
\begin{theorem}\label{T:2.3}
Let the pole point be $P$ of the one-parameter homothetic
hyperbolic motion $ {\mathbb{H} \mathord{\left/
 {\vphantom {a b}} \right.
 \kern-\nulldelimiterspace} \mathbb{H'}}$. During this motion the acceleration pole point
 $Q = \left( {q_1 ,q_2 } \right) \in \mathbb{H}$ is the hyperbolic number as
\begin{equation}\label{2.19}
{\bf q} = {\bf p} + \frac{{{\bf \dot p}\left( {\dot h + jh\dot
\varphi } \right)}}{{\ddot h + h\dot \varphi ^2  + j\left( {2\dot
h\dot \varphi  + h\ddot \varphi } \right)}}
\end{equation}
where $ \ddot h + h\dot \varphi ^2  \ne  \mp \left( {2\dot h\dot
\varphi  + h\ddot \varphi } \right)$.
\end{theorem}
\textbf{Special Case} In the special case of $h=1$, given the
condition that $\ddot \varphi ^2  - \dot \varphi ^4  \ne 0$, the
acceleration pole point of one-parameter hyperbolic motion is
\[
{\bf q} = {\bf p} + \frac{{{\bf \dot p}\left( {\dot \varphi \ddot
\varphi  - j\dot \varphi ^3 } \right)}}{{\ddot \varphi ^2  - \dot
\varphi ^4 }}
\]
which was given in \cite{Yuc}.

\section{Canonical Relative System for Homothetic Motion in the Hyperbolic Plane}

Let's consider an $\mathbb{A}$ plane which moves with regard to
$\mathbb{H}$ and $\mathbb{H'}$ hyperbolic planes, first one moving
and the second one fixed. Let's examine the motion of the
coordinate system $\left\{ {B;\;{\bf a}_{\bf 1} ,\;{\bf a}_{\bf 2}
} \right\}$ which defines hyperbolic plane $\mathbb{A}$, and
hyperbolic planes $\mathbb{H}$ and $\mathbb{H'}$ with regard to
the coordinate systems $\left\{ {O;\;{\bf h}_{\bf 1} ,\;{\bf
h}_{\bf 2} } \right\}$ and $\left\{ {O';\;{\bf h'}_1 ,\;{\bf h'}_2
} \right\}$. [See Figure 4. and 5.] If the vector $\overrightarrow
{OB}$ is defined by the hyperbolic number ${\bf b} = b_1  + jb_2$,
by applying the hyperbolic inner product, $b_1^2  - b_2^2  > 0$ or
$b_1^2  - b_2^2 < 0$ can be obtained. As seen in Figure 3.1. and
Figure 3.2. respectively, the vector $\overrightarrow {OB}$ can be
on the plane H-I or H-II in hyperbolic motion.
\hfil\scalebox{1}{\includegraphics{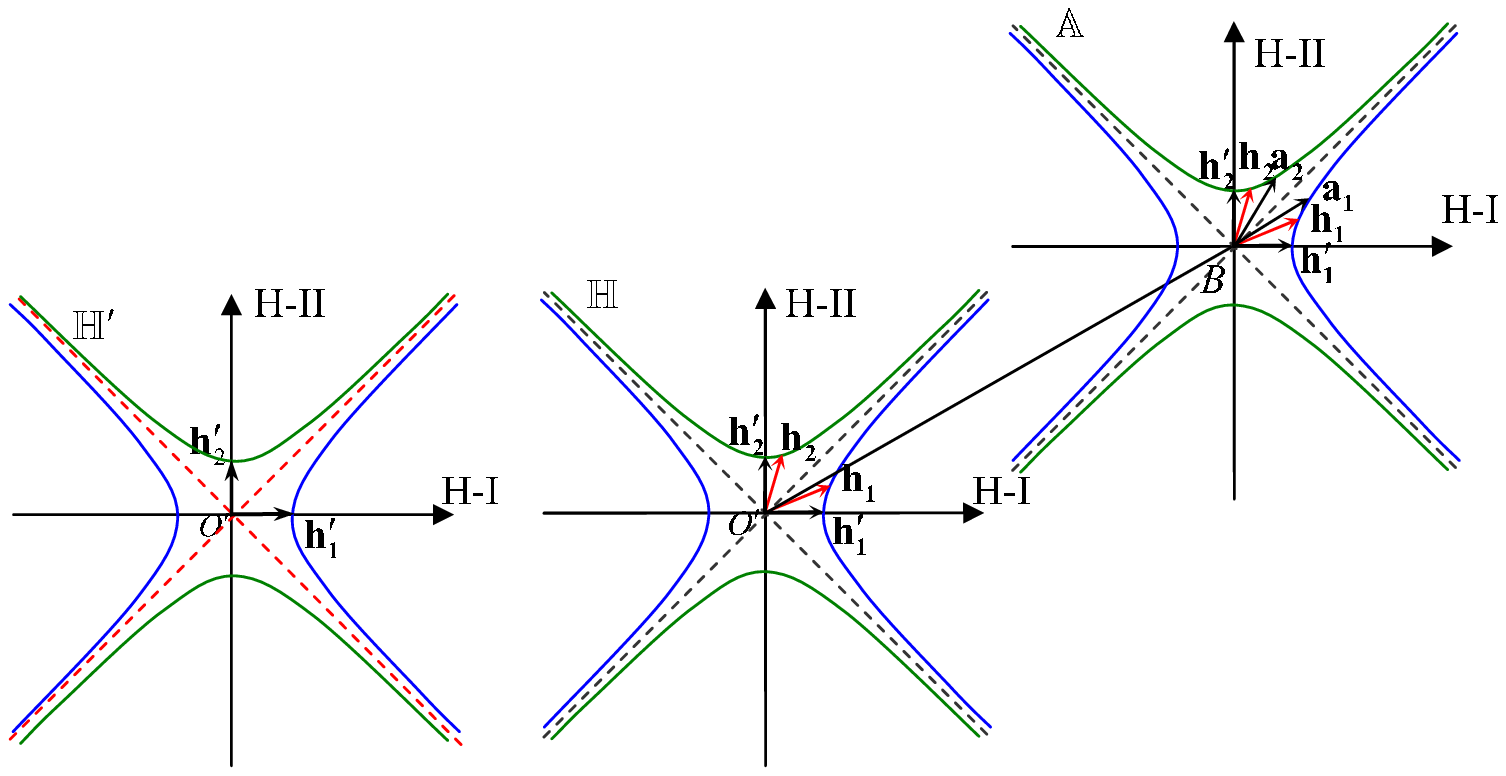}}\hfil
\begin{center}
\scriptsize{Figure  4. $\overrightarrow {OB}$ vector is on H-I
plane}
\end{center}
\hfil\scalebox{1}{\includegraphics{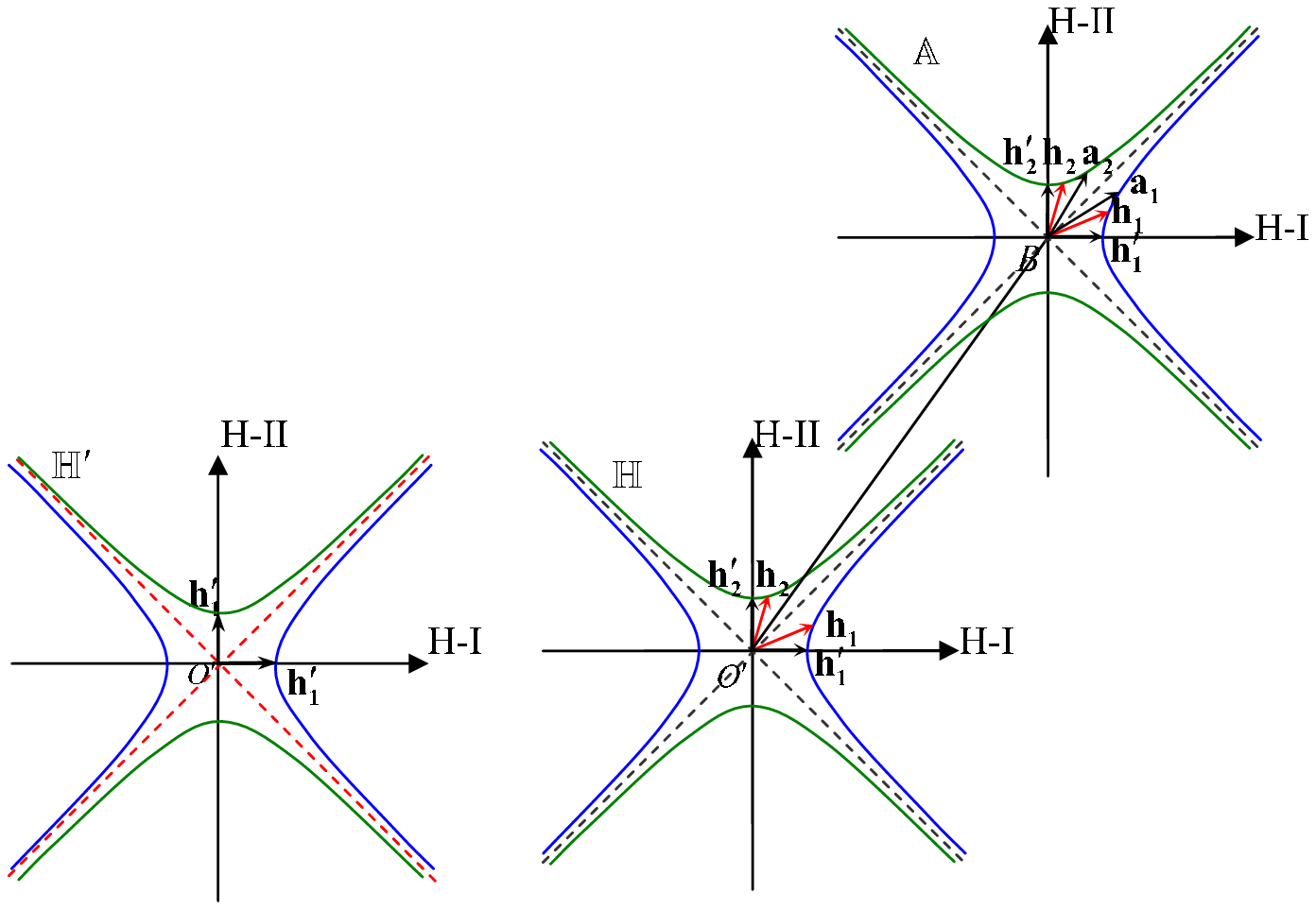}}\hfil
\begin{center}
\scriptsize{Figure 5. $\overrightarrow {OB}$ vector is on H-II
plane}
\end{center}
\normalsize The rotation angles of the one-parameter planar
hyperbolic motion ${\mathbb{\mathbb{A}} \mathord{\left/
 {\vphantom {A H}} \right.
 \kern-\nulldelimiterspace} \mathbb{H}}$ and  ${\mathbb{A} \mathord{\left/
 {\vphantom {A H}} \right.
 \kern-\nulldelimiterspace} \mathbb{H'}}$ are $\varphi$ and $\psi$, respectively. If the origin points of $O,\;B $ and $O',\;B $
 are coincident, then there exists following relations;
\[
\begin{array}{l}
 {\bf a}_{\bf 1}  = \cosh \varphi {\bf h}_{\bf 1}  + \sinh \varphi {\bf h}_{\bf 2}  \\
 {\bf a}_{\bf 2}  = \sinh \varphi {\bf h}_{\bf 1}  + \cosh \varphi {\bf h}_{\bf 2}  \\
 \end{array}
\]
and
\[
\begin{array}{l}
 {\bf a}_{\bf 1}  = \cosh \psi {\bf h'}_{\bf 1}  + \sinh \psi {\bf h'}_{\bf 2}  \\
 {\bf a}_{\bf 2}  = \sinh \psi {\bf h'}_{\bf 1}  + \cosh \psi {\bf h'}_{\bf 2}  \\
 \end{array}
\]
respectively [See Figure 6.].\\

\hfil\scalebox{1}{\includegraphics{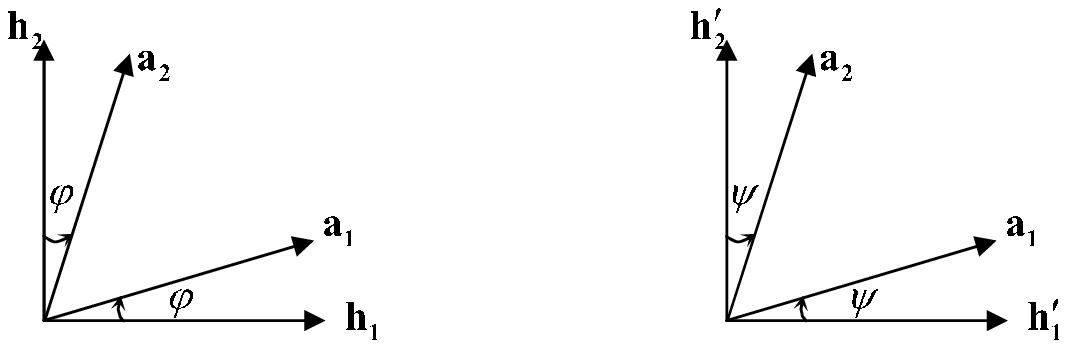}}\hfil\\
\begin{center}
\scriptsize{Figure 6.}
\end{center}
\normalsize
Let $X$  be a point with coordinates $ x_1 ,\;\,x_2$ on
the moving plane $\mathbb{A}$. If we denote the vectors
$\overrightarrow {BX}$, $\overrightarrow {OB}$ and
$\overrightarrow {OB'}$ with the hyperbolic numbers  $\tilde X =
x_1  + jx_2 $, ${\bf b} = b_1 + jb_2$ and ${\bf b'}=b'_1  +
jb'_2$, respectively; then we can write
\begin{equation}\label{3.1}
{\bf x} = \left( {{\bf b} + h{\bf \tilde x}} \right)e^{j\varphi }
\end{equation}
and
\begin{equation}\label{3.2}
{\bf x'} = \left( {{\bf b'} + h{\bf \tilde x'}} \right)e^{j\varphi
}
\end{equation}
where, $h = h\left( t \right) \ne {\rm constant}$ is the
homothetic scale of the motion and the hyperbolic numbers $\bf x$
and $\bf x'$ denote the point $X$ with respect to the coordinate
systems of $\mathbb{H}$ and $\mathbb{H'}$,
respectively.\\
The velocities of the motion with the help of the differentiation
of the equations (\ref{3.1}) and (\ref{3.2}) can be found. By
differentiating the equation (\ref{3.1})
\begin{equation}\label{3.3}
d{\bf x} = \left( {\sigma  + \left( {dh + jh\tau } \right){\bf
\tilde x} + hd{\bf \tilde x}} \right)e^{j\varphi }
\end{equation}
is obtained, in which
\begin{equation}\label{3.4}
\sigma  = \sigma _1  + j\sigma _2  = d{\bf b} + j{\bf b}d\varphi
\quad ,\quad \tau  = d\varphi
\end{equation}
and the relative velocity vector of $X$ (with respect to
$\mathbb{H}$) is $ {\bf V}_{\bf r}  = \frac{{d{\bf x}}}{{dt}}$.\\
If we assume the differentiation of the equation (\ref{3.2}),
\begin{equation}\label{3.5}
d'{\bf x} = \left( {\sigma ' + \left( {dh + jh\tau '} \right){\bf
\tilde x} + d{\bf \tilde x}} \right)e^{j\psi }
\end{equation}
can be obtained along with the equation
\begin{equation}\label{3.6}
\sigma ' = \sigma '_1  + j\sigma '_2  = d'{\bf b} + j{\bf b'}d\psi
\quad ,\quad \tau ' = d\psi .
\end{equation}
Also, the absolute velocity vector, that is, the velocity vector
of $X$ with respect to $\mathbb{H'}$, is ${\bf V_a}  = \frac{{d'{\bf x}}}{{dt}}$.\\
Here, $ \sigma _i ,\;\;\sigma '_i ,\;\;\left( {i = 1,2}
\right),\;\;\tau ,\;\;\tau ' $ are linear differential forms of
$t$ and are called Lorentzian Pfaffian forms of one-parameter
homothetic hyperbolic motion. The real parameter $t$ represents time.\\
If ${\bf V_r}= {\bf 0}$ and ${\bf V_a}= {\bf 0}$, the point $X$ is
fixed on the hyperbolic planes $\mathbb{H}$ and $\mathbb{H'}$,
respectively. Thus, the conditions of $X$ being fixed on the
$\mathbb{H}$ and $\mathbb{H'}$ planes are
\begin{equation}\label{3.7}
d{\bf \tilde x} =  - \frac{1}{h}\left( {\sigma  + \left( {dh +
jh\tau } \right){\bf \tilde x}} \right)
\end{equation}
and
\begin{equation}\label{3.8}
d{\bf \tilde x} =  - \frac{1}{h}\left( {\sigma ' + \left( {dh +
jh\tau '} \right){\bf \tilde x}} \right)
\end{equation}
respectively. If the equation (\ref{3.7}) is substituted into
equation (\ref{3.5}),
\begin{equation}\label{3.9}
d_f {\bf x} = \left[ {\left( {\sigma ' - \sigma } \right) +
jh\left( {\tau ' - \tau } \right){\bf \tilde x}} \right]e^{j\psi }
\end{equation}
can be obtained, where the sliding velocity vector of the point
$X$ is  ${\bf V_f } = \frac{{d_f {\bf x}}}{{dt}}$. Thus, following
can be easily obtained:
\begin{equation}\label{3.10}
d'{\bf x} = d_f {\bf x} + d{\bf x}
\end{equation}
This satisfies the relation between velocities which is given in
equation (\ref{2.9}).\\
Just to avoid translation, it is assumed that $\dot \varphi  \ne
0$ and  $ \dot \psi  \ne 0$. The rotation pole of the motion $
{\mathbb{H} \mathord{\left/
 {\vphantom {H {H'}}} \right.
 \kern-\nulldelimiterspace} {\mathbb{H'}}}$ is characterized by the sliding velocity  $P$ being zero. For that
reason, if $d_f {\bf x}={\bf 0}$, from the equation (\ref{3.9}),
the pole point $P$ of the one-parameter planar hyperbolic
homothetic motion is obtained as
\begin{equation}\label{3.11}
{\bf p}= j\frac{{\sigma ' - \sigma }}{{\tau  - \tau '}}
\end{equation}
and if Lorentzian coordinates are preferred on the condition that
$\overrightarrow {BP}  = {\bf p} = p_1  + jp_2 $, it can be
written
\begin{equation}\label{3.12}
p_1  = \frac{{\sigma '_2  - \sigma _2 }}{{\tau  - \tau '}}\quad
,\quad p_2  = \frac{{\sigma '_1  - \sigma _1 }}{{\tau  - \tau '}}
\end{equation}
which is given in \cite{Erg}.\\
In the ${\mathbb{H} \mathord{\left/
 {\vphantom {H {H'}}} \right.
 \kern-\nulldelimiterspace} {\mathbb{H'}}}$ one-parameter planar hyperbolic homothetic motion, moving and fixed
pole curves determine the geometric locus of the point $P$ in
$\mathbb{H}$ and $\mathbb{H'}$ planes, respectively. In other
words; $\left( P \right)$ and $\left( P'\right)$ are the
representation of the moving and fixed pole curves, respectively.
Also, the pole tangents can be
either on the plane H-I or H-II [See Figure 7.].\\

\hfil\scalebox{1}{\includegraphics{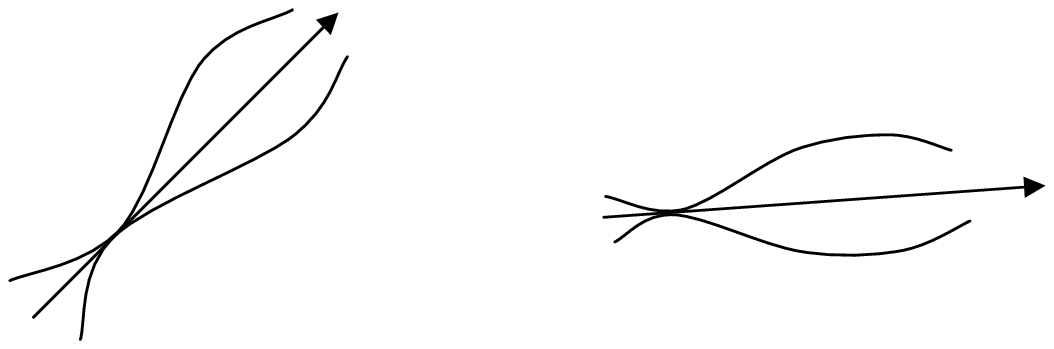}}\hfil\\
\begin{center}
\scriptsize{Figure 7.}
\end{center}
\normalsize
Let's first choose the pole tangents of the pole
curves $\left( P \right)$ and $\left( P'\right)$ on the plane H-II
because the same results would be obtained by following similar
operations on the
plane H-I.\\
\subsection{The Euler-Savary Formula for One-Parameter Planar
Hyperbolic Homothetic Motion}

Let's choose the moving plane $\mathbb{A}$, represented by the
coordinate system $ \left\{ {B;{\bf a}_{\bf 1} ,{\bf a}_{\bf 2} }
\right\}$,
in such way to meet the following conditions:\\
\textbf{i)} The origin of the system $B$ coincides with the
instantaneous rotation pole $P$\\
\textbf{ii) }The axis $ \left\{ {B;{\bf a}_2 } \right\}$ is the
pole tangent, that is, it coincides with the common tangent of the
pole curves $\left( P \right)$ and $\left( P'\right)$ (on the
plane H-II) [See Figure 8.].\\

\hfil\scalebox{1}{\includegraphics{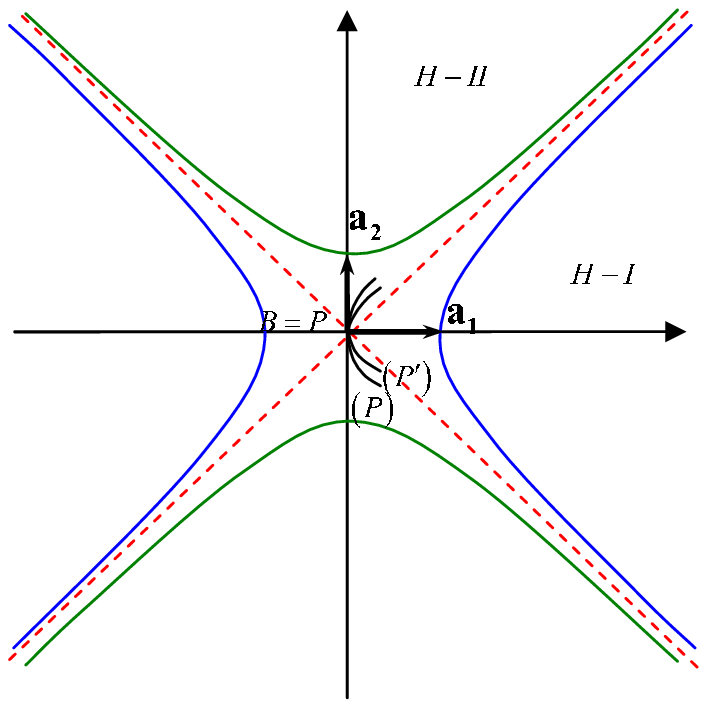}}\hfil
\begin{center}
\scriptsize{Figure 8.}
\end{center}
\normalsize When the condition (i) is considered: by using the
equation (\ref{3.12}),
\begin{equation}\label{3.13}
\sigma _1  = \sigma '_1 \quad ,\quad \sigma _2  = \sigma '_2
\end{equation}
are obtained. From the equations (\ref{3.4}) and  (\ref{3.6}),
\begin{equation}\label{3.14}
\begin{array}{*{20}c}
   {d{\bf b} = (d{\bf b} + j{\bf b}d\varphi )e^{j\varphi }  = \sigma e^{j\varphi } }  \\
   {d'{\bf b} = (d{\bf b'} + j{\bf b'}d\varphi )e^{j\psi }  = \sigma 'e^{j\psi } }  \\
\end{array}
\end{equation}
are found. If the equation (\ref{3.13}) and the last equation are
took into consideration:
\begin{equation}\label{3.15}
d{\bf p} = d'{\bf p} = d{\bf b} = d'{\bf b}
\end{equation}
is found.\\
Thus, the moving pole curve $\left( P \right)$, the pole tangent
of which is given, and the fixed pole curves $\left( P' \right)$
are rolling on each other without sliding.\\
The second condition, that is, the condition that the pole tangent
coincides with  ${\bf a}_{\bf 2}$, requires the coefficient of
${\bf a}_{\bf 1}$ to be zero. Here, $\sigma _1  = \sigma '_1  = 0$
and $\sigma  = j\sigma _2  = j\sigma '_2$ can be written.
Consequently, the derivative equations of the canonical relative
system $\{ P;{\bf a}_{\bf 1} ,{\bf a}_{\bf 2} \}$ are
\begin{equation}\label{3.16}
d{\bf a}_{\bf 1}  = \tau \,{\bf a}_{\bf 2}  = j\tau {\kern 1pt}
e^{j\varphi } \;,\quad d{\bf a}_{\bf 2}  = \tau \,{\bf a}_{\bf 1}
=  - {\kern 1pt} \tau {\kern 1pt} e^{j\varphi } \;,\quad d{\bf p}
= j\sigma _2 {\kern 1pt} {\bf a}_{\bf 1}  = \sigma {\kern 1pt}
e^{j\varphi }
\end{equation}
and
\begin{equation}\label{3.17}
d'{\bf a}_{\bf 1}  = \tau '{\bf a}_{\bf 2}  = j\tau 'e^{j\psi }
,\quad \quad d'{\bf a}_{\bf 2}  = \tau '{\bf a}_{\bf 1}  = j\tau
'e^{j\psi } ,\quad \quad d'{\bf p} = j\sigma '_2 {\kern 1pt} {\bf
a}_{\bf 1}  = \sigma {\kern 1pt} e^{j\psi } .
\end{equation}
Here $\sigma  = ds$ is the scalar arc element of the pole curves
$\left( P \right)$ and $\left( P'\right)$. $\tau$ is the
hyperbolic cotangent angle, that is, two neighboring tangent
angles of $\left( P \right)$. Thus, the curvature of $\left( P
\right)$ on the point $P$ is represented by $\frac{\tau }{\sigma }
= \frac{{d\varphi }}{{ds}}$.\\
Similarly, the curvature of the fixed pole curve $\left( P'
\right)$ on the point $P$ is $\frac{{\tau '}}{\sigma } =
\frac{{d\psi }}{{ds}}$ where $\tau'$ is the hyperbolic cotangent
angle.\\
The inverse values of these ratios
\begin{equation}\label{3.18}
r = \frac{\sigma }{\tau }
\end{equation}
and
\begin{equation}\label{3.19}
r' = \frac{\sigma }{{\tau '}}
\end{equation}
give the curvature radius of the pole curves $\left( P \right)$
and $\left( P'\right)$, respectively.\\
When $d\nu  = \tau ' - \tau$ is the infinitesimal small hyperbolic
instantaneous rotation angle, the moving hyperbolic plane
$\mathbb{\mathbb{H}}$, with respect to the fixed plane
$\mathbb{H'}$, rotates around the rotation pole $P$ as much as
this hyperbolic angle in the $dt$ time scale. Thus, the hyperbolic
angular velocity of the rotational motion of $\mathbb{H}$ with
respect to $\mathbb{H'}$ is
\begin{equation}\label{3.20}
\frac{{\tau ' - \tau }}{{dt}} = \frac{{d\nu }}{{dt}} = \mathop \nu
\limits^ \bullet.
\end{equation}
From the equations (\ref{3.18}), (\ref{3.19}), and the last
equation, the following can be written:
\begin{equation}\label{3.21}
\frac{{\tau ' - \tau }}{{dt}} = \frac{{d\nu }}{{dt}} =
\frac{1}{{r'}} - \frac{1}{r}.
\end{equation}
Let the direction of the unit tangent vector ${\bf a}_{\bf 2} $ be
in the direction determined by time-based pole curves $\left( P
\right)$ and $\left( P'\right)$. Let's choose the vector ${\bf
a}_{\bf 2} $ in such way to ensure that $\frac{{ds}}{{dt}} > 0$.
In this case, $r > 0$ as the curvature center of the moving pole
$\left( P \right)$ curve is at the right side of the directed pole
tangent $\left\{ {P;{\bf a}_{\bf 2} } \right\}$. Similarly, $r' >
0$.\\
According to the canonical relative system, the differentiation
${\bf x}$- the coordinates of which are $x_1 ,\;x_2$- with respect
to the planes $\mathbb{H}$ and $\mathbb{H'}$ are
\begin{equation}\label{3.22}
d{\bf x} = (\sigma  + \left( {dh + jh\tau } \right){\bf \tilde x}
+ hd{\bf \tilde x})e^{j\varphi }
\end{equation}
and
\begin{equation}\label{3.23}
d'{\bf x} = (\sigma ' + \left( {dh + jh\tau '{\bf \tilde x}}
\right) + hd{\bf \tilde x})e^{j\psi }
\end{equation}
respectively. If
\begin{equation}\label{3.24}
hd{\bf \tilde x} =  - \sigma  - \left( {dh + jh\tau } \right){\bf
\tilde x}
\end{equation}
then the point $X$ is fixed on the hyperbolic plane $\mathbb{H}$.
Similarly, if
\begin{equation}\label{3.25}
hd{\bf \tilde x} =  - \sigma  - \left( {dh + jh\tau '} \right){\bf
\tilde x}
\end{equation}
then the point $X$ is fixed on the hyperbolic plane $\mathbb{H'}$.
Also, the sliding velocity $\bf{V_f}$ of the movement ${\mathbb{H}
\mathord{\left/
 {\vphantom {H {H'}}} \right.
 \kern-\nulldelimiterspace} {\mathbb{H'}}}$ corresponds to the differentiation
\begin{equation}\label{3.26}
d_f {\bf x} = jh\left( {\tau ' - \tau } \right){\bf \tilde
x}e^{j\varphi } .
\end{equation}
Now, let's examine the curvature centers of the trajectory curves
drawn on their fixed plane by the points of moving planes in the
motion of ${\mathbb{H} \mathord{\left/ {\vphantom {H {H'}}}
\right.
 \kern-\nulldelimiterspace} {\mathbb{H'}}}$. In the canonical relative system, the points $X$, $X'$ having
the coordinates $x_1$, $x_2$ and $x'_1$, $x'_2$, respectively, are
situated, -together with the instantaneous rotation pole $P$ in
every $t$ moment on the instantaneous trajectory normal, which
belongs to $X$. Moreover, this curvature center can be considered
as the limit of the meeting point of the normals of the two
neighboring points on the curve. Thus,
\begin{equation}\label{3.27}
\begin{array}{l}
 \overrightarrow {PX}  = x_1  + jx_2  = \bf{x} \\
 \overrightarrow {PX'}  = x'_1  + jx'_2  = \bf{x'} \\
 \end{array}
\end{equation}
vectors have the same direction which passes through $P$. Then,
for the points $X$ and $X'$, the equation is
\begin{equation}\label{3.28}
\frac{\bf{x}}{{\bf{x'}}} = \frac{{x_1  + jx_2 }}{{x'_1  + jx'_2 }}
= \lambda \in\mathbb{R}.
\end{equation}
If the differential of this last equation is taken, then we get
\begin{equation}\label{3.29}
d{\bf xx'} - {\bf x}d{\bf x'} = 0.
\end{equation}
If the conditions that the point $X$ be fixed on the plane
$\mathbb{H}$ and the point $X'$ be fixed on the plane
$\mathbb{H'}$ are provided, then
\begin{equation}\label{3.30}
\sigma \left[ {{\bf x} - {\bf x'}} \right] + jh{\bf xx'}(\tau ' -
\tau ) = 0
\end{equation}
can be obtained. As the vectors $ \overrightarrow {PX}$,
$\overrightarrow {PX'}$ are on the plane H-II ,
\begin{equation}\label{3.31}
{\bf{x}} = aje^{j\alpha }
\end{equation}
and
\begin{equation}\label{3.32}
{\bf{x'}} = a'je^{j\alpha }.
\end{equation}
That is, $a$ and $a'$, respectively, represent the distance of the
points $X$ and $X'$ on the plane H-II from the rotation pole $P$.
Also, the angle $\alpha$ is bounded by the pole curves
$\overrightarrow {PX}  = \overrightarrow {PX'}$, [See Figure 9.]

\hfil\scalebox{1}{\includegraphics{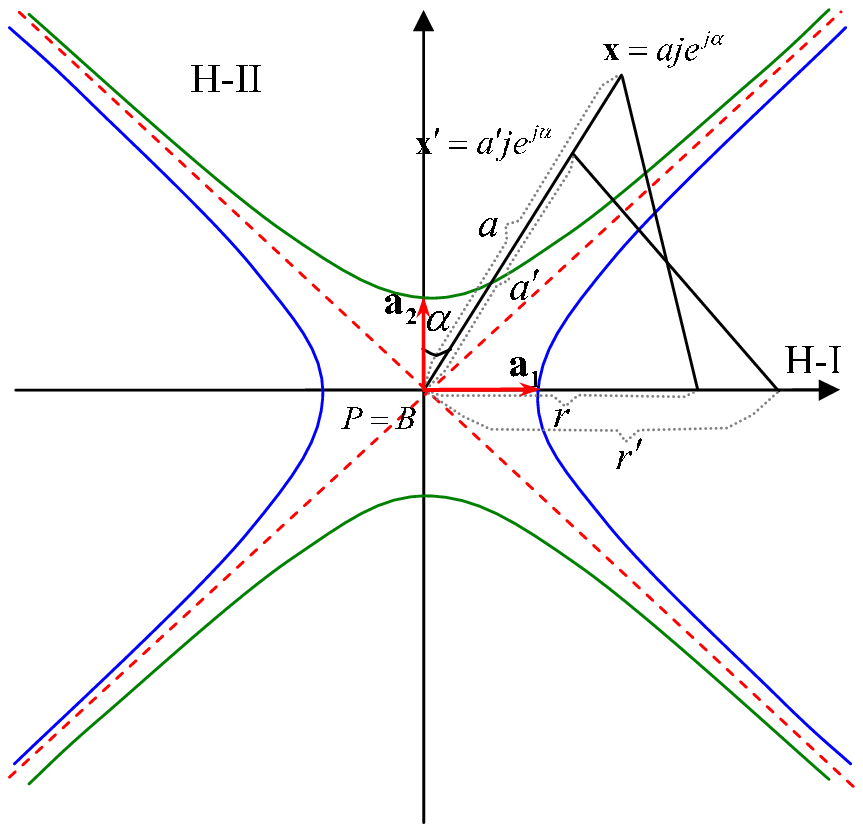}}\hfil
\begin{center}
\scriptsize{Figure 9.}
\end{center}
If the equations (\ref{3.31}) and (\ref{3.32}) are substituted
into equation (\ref{3.30}), then
\begin{equation}\label{3.33}
j \sigma (a - a') + jhaa'e^{j\alpha } (\tau ' - \tau ) = 0
\end{equation}
can be obtained, and if the equation (\ref{3.21}) is considered
together with this last equation,
\begin{equation}\label{3.34}
\left( {\frac{1}{a} - \frac{1}{{a'}}} \right)e^{ - j\alpha }  =
h\left( {\frac{1}{{r'}} - \frac{1}{r}} \right) = h\frac{{d\nu
}}{{ds}}.
\end{equation}
is found. Here, $r$ and $r'$ are the radii of curvature of the
pole curves $P$ and $P'$, respectively. $ds$ represents the scalar
arc element and $d\nu$ represents the infinitesimal hyperbolic
angle of the motion of the pole curves.\\
The equation (\ref{3.34}) is called the Euler-Savary formula for
one-parameter plane hyperbolic homothetic motion.\\
Consequently, the following theorem can be given.
\begin{theorem}\label{T:3.1}
Let $\mathbb{H}$ and $\mathbb{H'}$ be the moving and fixed
hyperbolic planes, respectively. A point $X$, assumed on
$\mathbb{H}$, draws a trajectory whose instantaneous center of
curvature is $X'$ on the plane $\mathbb{H'}$ in one-parameter
planar homothetic motion ${\mathbb{H} \mathord{\left/
 {\vphantom {H {H'}}} \right.
 \kern-\nulldelimiterspace} {\mathbb{H'}}}$. In the inverse homothetic motion of ${\mathbb{H} \mathord{\left/
 {\vphantom {H {H'}}} \right.
 \kern-\nulldelimiterspace} {\mathbb{H'}}}$, a point $X'$ assumed on $\mathbb{H'}$ draws a
trajectory whose center of curvature is $X$ on the plane
$\mathbb{H}$. The relation between the points $X$ and $X'$ is
given by the Euler-Savary formula given in the equation
(\ref{3.34}).
\end{theorem}
\textbf{Remark} Let's choose the moving plane $\mathbb{A}$
represented by the coordinate system $ \left\{ {B;{\bf a}_{\bf 1}
,{\bf a}_{\bf 2} } \right\}$
in such way to meet following conditions:\\
\textbf{i)} The origin of the system $B$ and the instantaneous
rotation pole $P$ coincide with each other, i.e. $B=P$, [See
Figure 10.]\\
\textbf{ii)} The axis $ \left\{ {B;{\bf a}_1 } \right\}$  is the
pole tangent, that is, it coincides with the common tangent of the
pole curves $\left( P \right)$ and $\left( P'\right)$ (on the
plane H-I)

\hfil\scalebox{1}{\includegraphics{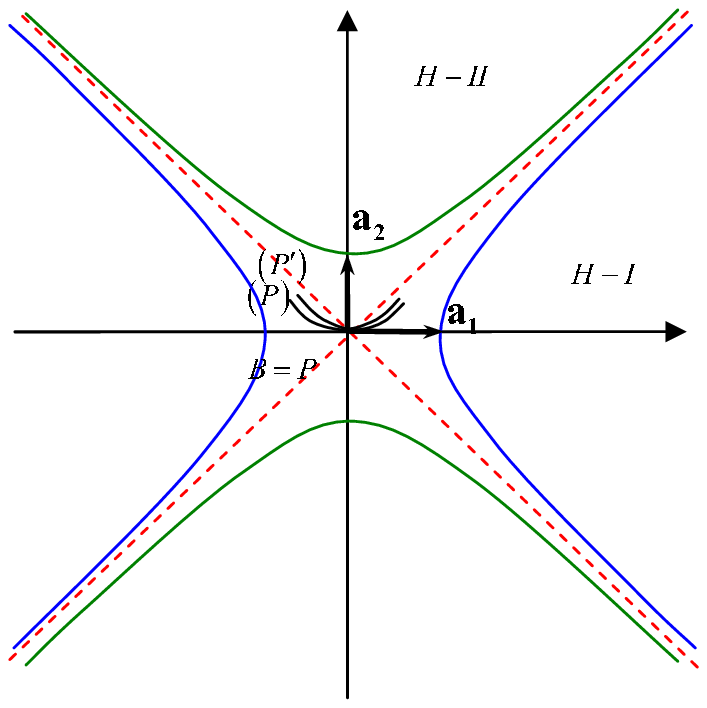}}\hfil
\begin{center}
\scriptsize{Figure 10.}
\end{center}
Thus, if the operations in section3.1. are performed considering
the conditions i) and ii), the Euler-Savary formula for
one-parameter planar hyperbolic homothetic motion remains
unchanged, that is, it is the same as in the equation
(\ref{3.34})[See Figure 11.]

\hfil\scalebox{1}{\includegraphics{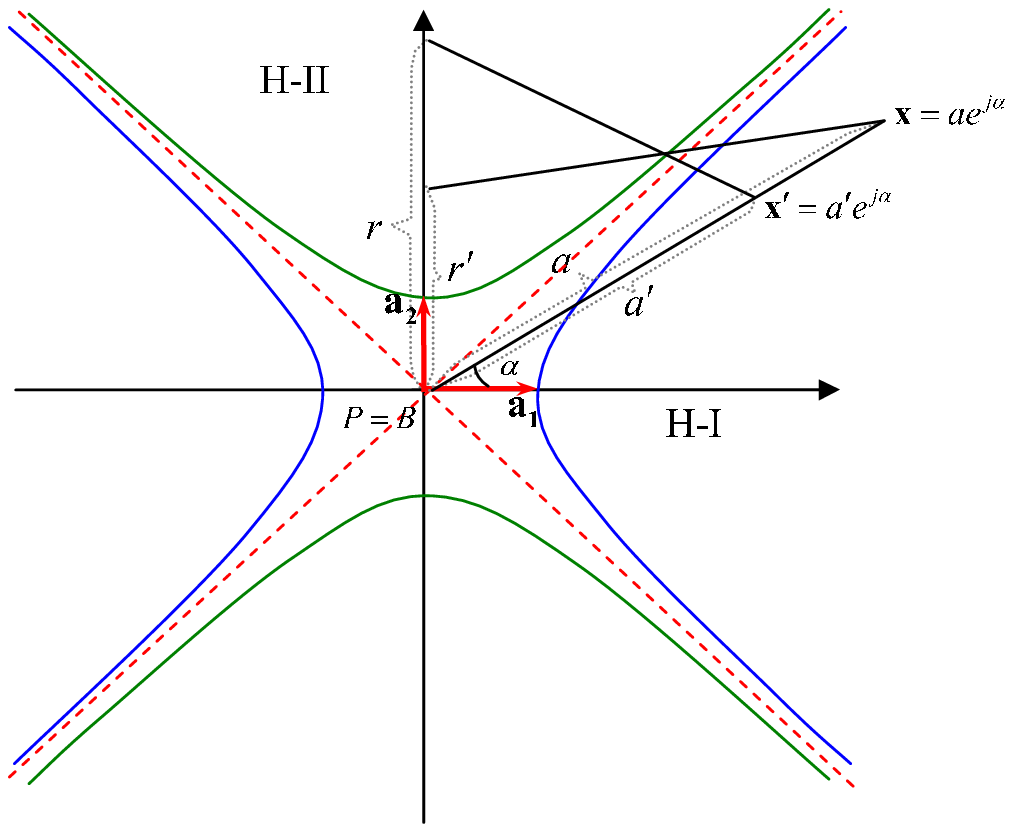}}\hfil
\begin{center}
\scriptsize{Figure 11.}
\end{center}

\end {document}